\documentclass[12pt]{article}
\usepackage{amsthm,amsfonts,amssymb,amscd,amsmath}

\begin{document}

\begin{center}
\Large \bf On the locus of non-rigid hypersurfaces
\end{center}
\vspace{1cm}

\centerline{Thomas Eckl and Aleksandr Pukhlikov}

\parshape=1
3cm 10cm \noindent {\small \quad \quad \quad
\quad\quad\quad\quad\quad\quad\quad {\bf }\newline We show that
the Zariski closure of the set of hypersurfaces of degree $M$ in
${\mathbb P}^{M}$, where $M\geq 5$, which are either not factorial
or not birationally superrigid, is of codimension at least $\binom{M-3}{2}+1$ in
the parameter space.

Bibliography: 21 titles.} \vspace{1cm}

\noindent {\bf 1. Formulation of the main result and scheme of the proof.}
Let ${\mathbb P}^{M}$, where $M\geq 5$, be the complex projective
space, ${\cal F}={\mathbb P}(H^0({\mathbb P}^{M},{\cal
O}_{{\mathbb P}^{M}}(M)))$ the space parametrizing hypersurfaces
of degree $M$. There are Zariski open subsets ${\cal F}_{\rm
reg}\subset {\cal F}_{\rm sm}\subset {\cal F}$, consisting of
hypersurfaces, regular in the sense of \cite{Pukh98b}, and smooth,
respectively. The well known theorem proven in \cite{Pukh98b}
claims that every regular hypersurface $V\in{\cal F}_{\rm reg}$ is
birationally superrigid. Let ${\cal F}_{\rm srigid}\subset{\cal
F}$ be the set of (possibly singular) hypersurfaces that are factorial and
birationally superrigid. The aim of this note is to show the
following claim.

{\bf Theorem 1.} {\it The Zariski closure $\overline{{\cal
F}\backslash{\cal F}_{\rm srigid}}$ of the complement is of
codimension at least $\binom{M-3}{2}+1$ in ${\cal F}$}.

Note that we do not discuss the question of whether ${\cal F}_{\rm
srigid}$ is open or not.

We prove Theorem 1, directly constructing a set in
${\cal F}$, every point of which corresponds to a factorial and
birationally superrigid hypersurface, with the Zariski closure of its complement of
codimension at least $\binom{M-3}{2}+1$. More precisely, let ${\cal F}_{{\rm
qsing}\geq r}$ be the set of hypersurfaces, every point of which
is either smooth or a quadratic singularity of rank at least $r$.
We do {\it not} assume that singularities are isolated, but it is
obvious that for $V\in{\cal F}_{{\rm qsing}\geq r}$ the following
estimate holds:
$$
\mathop{\rm codim}\mathop{\rm Sing}V\geq r-1.
$$
In particular, by the famous Grothendieck theorem (\cite[XI.Cor.3.14]{Gro}, \cite{CL}) any
$V\in{\cal F}_{{\rm qsing}\geq 5}$ is a factorial variety,
therefore a Fano variety of index 1:
$$
\mathop{\rm Pic}V={\mathbb Z}K_V,\,\,K_V=-H,
$$
where $H\in\mathop{\rm Pic}V$ is the class of a hyperplane
section.

It is easy to see (Proposition 2) that $\mathop{\rm codim}({\cal
F}\backslash{\cal F}_{{\rm qsing}\geq 5})\geq \binom{M-3}{2}+1$.

Denote by ${\cal F}_{{\rm reg,\, qsing}\geq 5}\subset{\cal
F}_{{\rm qsing}\geq 5}$ the subset, consisting of such Fano
hypersurfaces $V\in{\cal F}$ that:

(1) at every smooth point the regularity condition of
\cite{Pukh98b} is satisfied;

(2) through every singular point there are only finitely many
lines on $V$.

We obtain Theorem 1 from the following two facts.

{\bf Theorem 2.} {\it The codimension of the complement of ${\cal
F}_{{\rm reg,\, qsing}\geq 5}$ in ${\cal F}$ is at least
$\binom{M-3}{2}+1$} if $M \geq 5$.

{\bf Theorem 3.} {\it Every hypersurface $V\in{\cal F}_{{\rm
reg,\, qsing}\geq 5}$ is birationally superrigid}.

{\bf Proof of Theorem 2} is straightforward and follows the
arguments of \cite{Pukh98b, Pukh01}; it is given in Section 2.

{\bf Proof of Theorem 3} starts in the usual way
\cite{Pukh98b,Pukh01,Pukh07b}: take a mobile linear system
$\Sigma\subset|nH|$ on a hypersurface $V\in{\cal F}_{{\rm reg, \,
qsing}\geq 5}$. Assume that for a generic $D\in\Sigma$ the pair
$(V,\frac{1}{n}D)$ is not canonical, that is, the system $\Sigma$
has a {\it maximal singularity} $E\subset V^+$, where
$\varphi\colon V^+\to V$ is a birational morphism, $V^+$ a smooth
projective variety, $E$ a $\varphi$-exceptional divisor and the
{\it Noether-Fano inequality}
$$
\mathop{\rm ord}\nolimits_E\varphi^*\Sigma> na(E)
$$
is satisfied (see \cite{Pukh07b} for definitions and details). We
need to get a contradiction, which would immediately imply
birational superrigidity and complete the proof of Theorem 3.

We proceed in the standard way.

Let $D_1,D_2\in\Sigma$ be generic divisors and $Z=(D_1\circ D_2)$
the self-intersection of the system $\Sigma$. Further, let
$B=\varphi(E)$ be the centre of the maximal singularity $E$. If
${\mathop{\rm codim}}_V B=2$, then
$$
{\mathop{\rm codim}}_B(B\cap\mathop{\rm Sing}V)\geq 2,
$$
so we can take any curve $C\subset B$, $C\cap\mathop{\rm
Sing}V=\emptyset$, and applying \cite[Sec.3]{Pukh98b}, conclude that
$$
\mathop{\rm mult}\nolimits_C\Sigma\leq n.
$$
As $\mathop{\rm mult}_B\Sigma>n$, we get a contradiction. So we
may assume that ${\mathop{\rm codim}}_V B\geq 3$.

{\bf Proposition 1 (the $4n^2$-inequality).} {\it The following
estimate holds:}
$$
\mathop{\rm mult}\nolimits_BZ>4n^2.
$$

If $B\not\subset\mathop{\rm Sing}V$, then the $4n^2$-inequality is
a well known fact going back to the paper on the quartic
three-fold \cite{IM}, so in this case no proof is needed, see
\cite[Ch.~2]{Pukh07b} for details. Therefore we assume that
$B\subset\mathop{\rm Sing}V$. In that case Proposition 1 is a
non-trivial new result, proved below in Sec. 3. The proof makes use
of the fact that the condition of having at most quadratic
singularities of rank $\geq r$ is stable with respect to blow ups,
in some a bit subtle way. That fact is shown in Sec. 4.

Now we complete the proof of Theorem 3, repeating word for word
the arguments of \cite{Pukh98b}. Namely, we choose an irreducible
component $Y$ of the effective cycle $Z$, satisfying the
inequality
$$
\frac {\mathop{\rm mult}_oY}{\mathop{\rm deg}Y}>\frac{4}{M},
$$
where $o\in B$ is a point of general position. Applying the
technique of hypertangent divisors in precisely the same way as it
is done in \cite{Pukh98b} (see also \cite[Ch.~3]{Pukh07b}), we
construct a curve $C\subset Y$, satisfying the inequality
$\mathop{\rm mult}_oC>\mathop{\rm deg}C$, which is impossible. It
is here that we need the regularity conditions. This contradiction
completes the proof of Theorem 3.

{\bf Remark 1.} (i) $4n^2$-inequality is not true for a quadratic
singularity of rank $\leq 4$: the non-degenerate quadratic point of
a three fold shows that $2n^2$ is the best we can achieve.

(ii) Birational superrigidity of Fano  hypersurfaces with
non-degenerate quadratic singularities was shown in
\cite{Pukh03b}. Birational (super)rigidity of Fano hypersurfaces
with isolated singular points of higher multiplicities $3\leq m\leq M-2$
was proved in \cite{Pukh02a}, but the argument is really hard.
These two results show that the estimate for the codimension of
the non-rigid locus could most probably be considerably
sharpened.

(iii) There are a few other papers where various classes of
singular Fano varieties were studied from the viewpoint of their
birational rigidity. The most popular object was three-dimensional
quartics \cite{Pukh89c,CoMe,Me04,Sh08}. Other families were
investigated in \cite{Ch04b,Ch07a}. A family of Fano varieties
(Fano double spaces of index one) with a higher dimensional
singular locus was recently proven to be birationally superrigid
in \cite{Mul}.

(iv) A recent preprint of de Fernex \cite{dF12} proves birational
superrigidity of a class of Fano hypersurfaces of degree $M$ in
$\mathbb{P}^M$ with not necessarily isolated singularities without
assuming regularity. But the dimension of the singularity locus is
bounded by $\frac{1}{2}M-4$, and no estimate of the codimension of
the complement of this class is given.\vspace{0.3cm}


\noindent {\bf 2. The estimates for the codimension.} Let us prove Theorem 2.

First we discuss the regularity conditions in more details. Let $x$ be a smooth point on a hypersurface $V$ of degree $M$ in $\mathbb{P}^M$. Choose homogeneous coordinates $(X_0 : \ldots : X_M)$ on $\mathbb{P}^M$ such that $x = (1:0:\ldots:0)$. Then $V \cap \{X_0 \neq 0\}$ is the vanishing locus of a polynomial
$$
q_1 + \cdots + q_M
$$
where each $q_i$ is a homogeneous polynomial of degree $i$ in $M$ variables $X_1, \ldots, X_M$. The regularity condition of \cite{Pukh98b} states that $q_1, \ldots, q_{M-1}$ is a regular sequence in $\mathbb{C}[x_1, \ldots, x_M]$. In particular,
$$
{\mathop{\rm codim}}_{\mathbb{A}^M}(\{q_1 = \ldots = q_{M-1} = 0\}) = 1.
$$
Since all the vanishing loci $\{q_i = 0\}$ are cones with vertex in $x$, the set $\{q_1 = \ldots = q_{M-1} = 0\}$ must consist of a finite number of lines through $x$. Hence there also is only a finite number of lines on $V$ through $x$.

If $x$ is a singular point on $V$ then $q_1 \equiv 0$. The regularity condition (2) is equivalent to
$$
{\mathop{\rm codim}}_{\mathbb{A}^M}(\{q_2 = \ldots = q_M = 0\}) = 1.
$$
since $\{q_2 = \ldots = q_M = 0\} \subset V$, and because of homogeneity every line through $x$ on $V$ also lies in $\{q_i = 0\}$.

It is not known whether the set ${\cal F}_{\rm reg}$ is Zariski-open in ${\cal F}$, but it certainly contains a Zariski-open subset of ${\cal F}$. The codimension in ${\cal F}$ of its complement ${\cal F}\setminus{{\cal F}_{\rm reg}}$ is defined as the codimension of the Zariski closure of the complement. On the other hand, ${\cal F}_{{\rm qsing}\geq 5}$ is certainly Zariski-open, hence ${\cal F}\setminus{{\cal F}_{{\rm qsing}\geq 5}}$ is Zariski-closed. We have
$$
{\mathop{\rm codim}}_{\cal F}({\cal F}\setminus{{\cal F}_{{\rm reg,\, qsing}\geq 5}}) = \min ({\mathop{\rm codim}}_{\cal F}({\cal F}\setminus{{\cal F}_{\rm reg}}), {\mathop{\rm codim}}_{\cal F}({\cal F}\setminus{{\cal F}_{{\rm qsing}\geq 5}})).
$$
Hence the estimate of Theorem 2 follows from the following two propositions:

{\bf Proposition 2.} {\it The codimension of the complement of ${\cal F}\setminus{{\cal F}_{{\rm qsing}\geq 5}}$ in ${\cal F}$ is at least $\binom{M-3}{2}+1$ if $M \geq 5$.}

{\bf Proposition 3.} {\it The codimension of the (Zariski closure of the) complement of ${\cal F}\setminus{{\cal F}_{\rm reg}}$ in ${\cal F}$ is at least $\frac{M(M-5)}{2}+4$ if $M \geq 5$.}
\vspace{0.3cm}

{\bf Proof of Proposition 2.} Let $S_M := \mathbb{P}^{\binom{M+1}{2} - 1}$ be the projectivized space of all symmetric $M \times M$-matrices with complex entries. Let $S_{M,r}$ be the projectivized algebraic subset of $M \times M$ symmetric matrices of rank $\leq r$. The locus $Q_r(P)$ of hypersurfaces $H \in {\cal F}$ with $P \in H$ a singularity that is at least a quadratic point of rank at most $r$ has codimension in ${\cal F}$ equal to
\begin{eqnarray*}
\mathrm{codim}_{\cal F} Q_r(P) & = & 1 + M + \mathrm{codim}_{S_M} S_{M,r} = 1 + M + \dim S_M - \dim S_{M,r} = \\
                               & = & M + \binom{M+1}{2} - \dim S_{M,r}.
\end{eqnarray*}
Let $G(M-r,M)$ be the Grassmann variety parametrizing $(M-r)$-dimensional subspaces of $\mathbb{C}^M$. To calculate $\dim S_{M,r}$ we consider the incidence correspondence (see~\cite[Ex.12.4]{Harr:AG})
\[ \Phi := \left\{ (A, \Lambda) : \Lambda^T \cdot A = A \cdot \Lambda = 0 \right\} \subset S_M \times G(M-r,M). \]
Since the fibers of the natural projection $\pi_2: \Phi \rightarrow G(M-r,M)$ is given by a linear subspace of $S_M$ of dimension $\binom{r+1}{2} - 1$, the variety $\Phi$ is irreducible of
\[ \dim \Phi = \binom{r+1}{2} - 1 + r(M-r). \]
Since on the other hand the natural projection $\pi_1: \Phi \rightarrow S_M$ is generically $1:1$ onto $S_{M,r}$, $\dim \Phi = \dim S_{M,r}$.

\noindent Consequently, since the $Q_r(P)$ cover $Q_r$ and $P$ varies in $\mathbb{P}^M$,
\[ \mathrm{codim}_{\cal F} Q_r \geq \mathrm{codim}_{\cal F} Q_r(P) - M = \binom{M-r+1}{2} + 1. \]
This completes the proof of Proposition 2.

\vspace{0.3cm}

\noindent For $r=4$, we have $\mathrm{codim}_{\cal F} Q_4 \geq M$ if
\[ \mathrm{codim}_{\cal F} Q_4 - M \geq \binom{M-3}{2} + 1 - M = \frac{(M-2)(M-7)}{2} \geq 0, \]
hence if $M \geq 7$.

\vspace{0.3cm}
{\bf Proof of Proposition 3.} Let $\Phi = \{(x,H):x \in H\} \subset \mathbb{P}^M \times {\cal F}$ be the incidence variety of hypersurfaces of degree $M$ in $\mathbb{P}^M$. Let $\Phi_{\rm reg}$ be the subset of pairs $(x,H)$ satisfying the regularity conditions. Note that the Zariski closure $\overline{\Phi\setminus\Phi_{\rm reg}}$ in $\Phi$ maps onto the Zariski closure $\overline{{\cal F}\setminus{{\cal F}_{\rm reg}}}$ in ${\cal F}$. The fiber of $\Phi_{\rm reg}$ over a point $x \in \mathbb{P}^M$ under the natural projection $\pi_1: \mathbb{P}^M \times {\cal F} \rightarrow \mathbb{P}^M$ can be described as
\[ \Phi_{\rm reg}(x) := \{H: x \in H\ {\rm satisfies\ the\ regularity\ conditions}\} \subset {\cal F}. \]
Choosing homogeneous coordinates $(X_0:\ldots:X_M)$ on $\mathbb{P}^M$ such that $x = \mbox{(1:0:\ldots:0)}$ we can write ${\cal F} = \mathbb{P}H^0(\mathbb{P}^M, {\cal O}_{\mathbb{P}^M}(M))$ as a projectivized product
\[ {\cal F} = \mathbb{P}(\bigoplus_{i=0}^M {\cal P}_{i,M}\cdot X_0^{M-i}), \]
where the ${\cal P}_{i,M}$ are the vector spaces of homogeneous polynomials in $X_1, \ldots, X_M$ of degree $i$. In particular, the $\pi_1$-fiber $\Phi(x)$ of $\Phi$ over $x$ is $\mathbb{P}(\bigoplus_{i=1}^M {\cal P}_{i,M}\cdot X_0^{M-i})$.

For another point $x^\prime \in \mathbb{P}^M$ also choose homogeneous coordinates $(X_0^\prime:\ldots:X_M^\prime)$ on $\mathbb{P}^M$ such that in these new coordinates $x^\prime = (1:0:\ldots:0)$. Then the projective-linear automorphism on $\mathbb{P}^M$ given by the coordinate change from $(X_0:\ldots:X_M)$ to $(X_0^\prime:\ldots:X_M^\prime)$ maps a polynomial $F(X_0, \ldots, X_M)$ to the polynomial $F(X_0^\prime,\ldots,X_M^\prime)$. In particular, the induced linear automorphism on the affine cone $H^0(\mathbb{P}^M, {\cal O}_{\mathbb{P}^M}(M))$ over ${\cal F}$ maps the product structure $\prod_{i=0}^M {\cal P}_{i,M}\cdot X_0^{M-i}$ onto the product structure
$\prod_{i=0}^M {\cal P}^\prime_{i,M}\cdot (X_0^\prime)^{M-i}$. Hence the induced projective-linear automorphism on ${\cal F}$ maps $\Phi(x)$ onto $\Phi(x^\prime)$ and $\Phi_{\rm reg}(x)$ to $\Phi_{\rm reg}(x^\prime)$ because the regularity conditions only depend on these product structures.

Consequently, the $\pi_1$-fibers of the Zariski closure $\overline{\Phi\setminus\Phi_{\rm reg}}$ are the Zariski closure $\overline{\Phi(x)\setminus\Phi_{\rm reg}(x)}$, hence
\[ \dim \overline{\Phi\setminus\Phi_{\rm reg}} = \dim \overline{\Phi(x)\setminus\Phi_{\rm reg}(x)} + M. \]
Since $\dim \overline{{\cal F}\setminus{\cal F}_{\rm reg}} \leq \dim \overline{\Phi\setminus\Phi_{\rm reg}}$ we conclude
\begin{eqnarray*}
{\mathop{\rm codim}}_{\cal F} \overline{{\cal F}\setminus{\cal F}_{\rm reg}} & \geq & \dim {\cal F} - \dim \overline{\Phi\setminus\Phi_{\rm reg}} = {\mathop{\rm codim}}_{\cal F} \overline{\Phi(x)\setminus\Phi_{\rm reg}(x)} - M = \\
 & = & {\mathop{\rm codim}}_{\Phi(x)} \overline{\Phi(x)\setminus\Phi_{\rm reg}(x)} - (M-1).
\end{eqnarray*}

Let $\widetilde{\Phi}(x) = \prod_{i=1}^M {\cal P}_{i,M}$ and $\widetilde{\Phi}_{\rm reg}(x)$ be the preimages of $\Phi(x)$, $\Phi_{\rm reg}(x)$ in the affine cone $H^0(\mathbb{P}^M, {\cal O}_{\mathbb{P}^M}(M)) = \prod_{i=0}^M {\cal P}_{i,M}$ over ${\cal F}$. Obviously,
\[ {\mathop{\rm codim}}_{\Phi(x)} \overline{\Phi(x)\setminus\Phi_{\rm reg}(x)} = {\mathop{\rm codim}}_{\widetilde{\Phi}(x)} \overline{\widetilde{\Phi}(x)\setminus{\widetilde{\Phi}}_{\rm reg}(x)}. \]
$\widetilde{\Phi}(x)\setminus{\widetilde{\Phi}}_{\rm reg}(x)$ consists of a subset $S_1$ Zariski-closed in ${\cal P}_{1,M}^\ast \times \prod_{i=2}^M {\cal P}_{i,M}$ of polynomials $q_1 + \cdots + q_M$ not satisfying regularity condition (1), where ${\cal P}_{1,M}^\ast = {\cal P}_{1,M}\setminus\{0\}$, and a Zariski-closed subset $S_2$ of $\{0\} \times \prod_{i=2}^M {\cal P}_{i,M}$ of polynomials $q_2 + \cdots + q_M$ not satisfying regularity condition (2). Hence $\overline{\widetilde{\Phi}(x)\setminus{\widetilde{\Phi}}_{\rm reg}(x)}$ is the union of the Zariski closure of $S_1$ in $\widetilde{\Phi}(x)$ and $S_2$. Consequently,
\[ {\mathop{\rm codim}}_{\widetilde{\Phi}(x)} \overline{\widetilde{\Phi}(x)\setminus{\widetilde{\Phi}}_{\rm reg}(x)} = \min ({\mathop{\rm codim}}_{{\cal P}_{1,M}^\ast \times \prod_{i=2}^M {\cal P}_{i,M}} S_1, {\mathop{\rm codim}}_{\widetilde{\Phi}(x)} S_2). \]
For $1 \leq j < i \leq M$ let $\pi_{i,j}: {\cal P}_{1,M}^\ast \times \prod_{k=2}^i {\cal P}_{k,M} \rightarrow {\cal P}_{1,M}^\ast \times \prod_{k=2}^j {\cal P}_{k,M}$ be the natural projection. Following the notations in \cite{Pukh98b} we set for $k=2, \ldots, M-1$
\[ Y_k := \{(q_1, \ldots, q_k) \in {\cal P}_{1,M}^\ast \times \prod_{i=2}^k {\cal P}_{i,M} : {\mathop{\rm codim}}_{\mathbb{P}^M} \{q_1 = \ldots = q_k = 0\} < k\}, \]
\[ R_k := ({\cal P}_{1,M}^\ast \times \prod_{i=2}^k {\cal P}_{i,M})\setminus Y_k, \]
\[ \mu_k := \min_{(q_1, \ldots, q_{k-1}) \in R_{k-1}} {\mathop{\rm codim}}_{\pi_{k,k-1}^{-1}(q_1, \ldots, q_{k-1})} (\pi_{k,k-1}^{-1}(q_1, \ldots, q_{k-1}) \cap Y_k). \]
$S_1$ can be stratified into disjoint subsets
\[ S_1 = \bigcup_{i=2}^{M-1} \pi_{M,i}^{-1}(Y_i) \cap \pi_{M,i-1}^{-1}(R_{i-1}). \]
Each stratum $\pi_{M,i}^{-1}(Y_i) \cap \pi_{M,i-1}^{-1}(R_{i-1})$ is Zariski-closed in $\pi_{M,i-1}^{-1}(R_{i-1})$. Hence
\begin{eqnarray*}
{\mathop{\rm codim}}_{{\cal P}_{1,M}^\ast \times \prod_{i=2}^M {\cal P}_{i,M}} S_1 & = & \min_{2 \leq i \leq M-1} {\mathop{\rm codim}}_{\pi_{M,i-1}^{-1}(R_{i-1})} (\pi_{M,i}^{-1}(Y_i) \cap \pi_{M,i-1}^{-1}(R_{i-1})) \\
 & \geq & \min_{2 \leq i \leq M-1} \mu_i.
\end{eqnarray*}

\noindent In the same way as for $S_1$ we obtain
\[ {\mathop{\rm codim}}_{\prod_{i=2}^M {\cal P}_{i,M}} S_2 \geq \min_{2 \leq i \leq M} \nu_i, \]
where
\[ \nu_k := \min_{(q_2, \ldots, q_{k-1}) \in Q_{k-1}} {\mathop{\rm codim}}_{\sigma_{k,k-1}^{-1}(q_2, \ldots, q_{k-1})} (\sigma_{k,k-1}^{-1}(q_2, \ldots, q_{k-1}) \cap Z_k), \]
\[ Q_k := \prod_{i=2}^k {\cal P}_{i,M}\setminus Z_k, \]
\[ Z_k := \{(q_2, \ldots, q_k) \in \prod_{i=2}^k {\cal P}_{i,M} : {\mathop{\rm codim}}_{\mathbb{P}^M} \{q_2 = \ldots = q_k = 0\} < k-1\} \]
and $\sigma_{k,k-1}: \prod_{i=2}^k {\cal P}_{i,M} \rightarrow \prod_{i=2}^{k-1} {\cal P}_{i,M}$ is the natural projection. Consequently,
\[ {\mathop{\rm codim}}_{\prod_{i=1}^M {\cal P}_{i,M}} S_2 \geq \min_{2 \leq i \leq M} \nu_i + M, \]
because $\dim {\cal P}_{1,M} = M$. Using the technique of \cite{Pukh98b},
\[ \mu_i \geq \binom{M}{i}, i = 2, \ldots, M-1,\ {\rm and\ } \nu_j \geq \binom{M+1}{j}, j = 2, \ldots, M. \]
Unfortunately these estimates are too weak for our purposes if $i = M-1$ and $j = M$. Using the technique of \cite{Pukh01} we obtain a better estimate for
\begin{eqnarray*}
\lefteqn{{\mathop{\rm codim}}_{\pi_{M,M-2}^{-1}(R_{M-2})} \pi_{M,M-2}^{-1}(R_{M-2}) \cap \pi_{M,M-1}^{-1}(Y_{M-1}) = \hspace{1.5cm}} \\
& & \hspace{1.5cm} {\mathop{\rm codim}}_{\pi_{M-1,M-2}^{-1}(R_{M-2})} \pi_{M-1,M-2}^{-1}(R_{M-2}) \cap Y_{M-1}.
\end{eqnarray*}
First of all, $\pi_{M-1,M-2}^{-1}(R_{M-2}) \cap Y_{M-1}$ fibers over ${\cal P}_{1,M}^\ast = R_1$, hence the codimension is at least the minimal codimension in a fiber. So we can fix a $q_1 \in R_1$ and choose affine coordinates $X_1, \ldots, X_M$ such that $q_1 = X_1$. Restricting the $q_2, \ldots, q_{M-1}$ to $\{X_1=0\} \cong \mathbb{A}^{M-1}$ we obtain homogeneous polynomials in the variables $X_2, \ldots, X_M$. Hence their vanishing sets can be projectivized in $\mathbb{P}^{M-2}$, and setting
\[ R_{M-3}^\prime := \{(q_2, \ldots, q_{M-2}): {\mathop{\rm codim}}_{\mathbb{P}^{M-2}}(\{q_2 = \ldots = q_{M-2} = 0\}) = M-3\} \subset \prod_{i=2}^{M-2} {\cal P}_{i,M-1}^\prime, \]
\[ Y_{M-2}^\prime := \{(q_2, \ldots, q_{M-1}): {\mathop{\rm codim}}_{\mathbb{P}^{M-2}}(\{q_2 = \ldots = q_{M-1} = 0\}) < M-2\} \subset \prod_{i=2}^{M-1} {\cal P}_{i,M-1}^\prime \]
we want to determine a lower bound for
\[ {\mathop{\rm codim}}_{(\pi_{M-1,M-2}^\prime)^{-1}(R_{M-3}^\prime)} (\pi_{M-1,M-2}^\prime)^{-1}(R_{M-3}^\prime) \cap Y_{M-2}^\prime. \]
Here, ${\cal P}_{i,M-1}^\prime$ is the space of homogeneous polynomials of degree $i$ in $M-1$ variables $X_2, \ldots, X_M$ and $\pi_{M-1,M-2}^\prime: \prod_{i=1}^{M-1} {\cal P}_{i,M-1}^\prime \rightarrow \prod_{i=1}^{M-2} {\cal P}_{i,M-1}^\prime$ is the natural projection.

For each tuple $(q_2, \ldots, q_{M-2}) \in R_{M-3}^\prime$, integers $2 \leq b \leq M-2$ and $2 \leq i_1 < \ldots < i_{b-1} \leq M-2$, there exists a $b$-dimensional linear subspace $L_b \subset \mathbb{P}^{M-2}$ such that $\{q_{i_1} = \ldots = q_{i_{b-1}} = 0\} \cap L_b \subset \mathbb{P}^{M-2}$ has only $1$-dimensional components. Vice versa, a tuple $(q_2, \ldots, q_{M-1})$ lies in $(\pi_{M-1,M-2}^\prime)^{-1}(R_{M-3}^\prime) \cap Y_{M-2}^\prime$ if for each $1$-dimensional irreducible component $B \subset \{q_2 = \ldots = q_{M-2} = 0\}$ spanning the linear subspace $\langle B \rangle \subset \mathbb{P}^{M-2}$ of dimension $b$ there exist integers $2 \leq i_1 < \ldots < i_{b-1} \leq M-2$ such that $\{q_{i_1} = \ldots = q_{i_{b-1}} = 0\} \cap \langle B \rangle$ contains $B$ as a $1$-dimensional irreducible component and $q_{i|B} \equiv 0$ for all $i \in \{2, \ldots, M-1\}\setminus\{i_1, \ldots, i_{b-1}\}$ (hence for all $2 \leq i \leq M-1$).

In the terminology of \cite{Pukh01} $q_{i_1}, \ldots, q_{i_{b-1}}$ is called a \textit{good sequence} for $B \subset \langle B \rangle$. Its existence can be shown inductively, using the regularity condition defining $R_{M-3}^\prime$.

If $b=1$, the $i_1, \ldots, i_{b-1}$ do not exist, and the condition restricts to
\[ q_{2|B} \equiv \ldots \equiv q_{M-1|B} \equiv 0 \]
on the line $B = \langle B \rangle$.

We can cover $(\pi_{M-1,M-2}^\prime)^{-1}(R_{M-3}^\prime) \cap Y_{M-2}^\prime$ by subsets $Z(b; i_1, \ldots, i_{b-1}; L_b)$ consisting of all tuples $(q_2, \ldots, q_{M-1}) \in (\pi_{M-1,M-2}^\prime)^{-1}(R_{M-3}^\prime)$ such that
\[ \dim \{q_{i_1} = \ldots = q_{i_{b-1}} = 0\} \cap L_b = 1, \]
$\{q_{i_1} = \ldots = q_{i_{b-1}} = 0\} \cap L_b$ contains irreducible components linearly spanning $L_b$ and $q_i \equiv 0$ on such a component, for each $2 \leq i \leq M-1$. Here, $1 \leq b \leq M-2$, $2 \leq i_1 < \ldots < i_{b-1} \leq M-2$, and the $L_b$ are paramatrized by the (projective) Grassmann variety $\mathbb{G}(b,M-2)$ of $b$-dimensional linear subspaces $L_b \subset \mathbb{P}^{M-2}$. For $b = 1$,
\[ Z(1;L) := \{(q_2, \ldots, q_{M-1}): q_{i|L} \equiv 0, 2 \leq i \leq M-1\}. \]
All these subsets are Zariski-closed in varying Zariski-open subsets of $(\pi_{M-1,M-2}^\prime)^{-1}(R_{M-3}^\prime)$.

For $b > 1$ they fiber surjectively onto $\prod_{k=1}^{b-1} {\cal P}_{i_k}^\prime$. Hence their codimension is estimated by a lower bound for each given $q_{i_1}, \ldots, q_{i_{b-1}}$, of the codimension  of all tuples of $q_i$, $i \in \{1, \ldots, M-1\}\setminus\{i_1, \ldots, i_{b-1}\}$, such that $q_{i|B} \equiv 0$ on an irreducible curve $B$ linearly spanning $L_b$. To find such a lower bound choose homogeneous coordinates $(X_2: \ldots : X_M)$ such that
\[ L_b = \{X_{b+3} = \ldots = X_M = 0\}. \]
Then $q_i \in {\cal P}_{i,M-1}^\prime$ cannot vanish on an irreducible curve $B$ linearly spanning all of $L_b$ if $q_{i|L_b}$ is of the form
\[ \prod_{k=1}^i (a_{k,2}X_2 + \cdots + a_{k,b+2}X_{b+2}). \]
Consequently the codimension of all $q_i \in {\cal P}_{i,M-1}^\prime$ vanishing on such a curve $B$ is at least the dimension of the space of polynomials in this form, that is $b \cdot i + 1$. Here, $b$ is the dimension of the space of hyperplanes in $\mathbb{P}^b$. It follows that the codimension of $Z(b; i_1, \ldots, i_{b-1}; L_b)$ in (a Zariski-open subset of) $(\pi_{M-1,M-2}^\prime)^{-1}(R_{M-3}^\prime)$ is at least
\begin{eqnarray*}
\sum_{\stackrel{2 \leq i \leq M-1}{i \neq i_1, \ldots, i_{b-1}}} (b \cdot i + 1) & \geq & b \cdot (2 + \cdots + (M-1-b) + (M-1)) + (M-1-b) \\
 & = & b \cdot \frac{(M-1-b)(M-b)}{2} + (b+1)(M-1) - 2b.
\end{eqnarray*}
Similarly, the codimension of $Z(1;L)$ in $(\pi_{M-1,M-2}^\prime)^{-1}(R_{M-3}^\prime)$ is at least
\[ 3 + \ldots + M = \frac{M(M+1)}{2} - 3 \]
because $i+1$ is the codimension of the set of polynomials $q_i \in {\cal P}_{i,M-1}^\prime$ vanishing on the line $L \subset \mathbb{P}^{M-2}$.

Taking all these data together
\[ {\mathop{\rm codim}}_{(\pi_{M-1,M-2}^\prime)^{-1}(R_{M-3}^\prime)} (\pi_{M-1,M-2}^\prime)^{-1}(R_{M-3}^\prime) \cap Y_{M-2}^\prime \] must be at least the minimum of the numbers
\begin{eqnarray*}
\lefteqn{b \cdot \frac{(M-1-b)(M-b)}{2} + (b+1)(M-1) - 2b - (b+1)(M-2-b)} \hspace{5cm} \\
 & = & b \cdot \frac{(M-1-b)(M-b)}{2} + b^2 + 1, 2 \leq b \leq M-2,
\end{eqnarray*}
and
\[ \frac{M(M+1)}{2} - 3 - 2(M-3) = \frac{M(M-3)}{2} + 3. \]
Here, $(b+1)(M-2-b)$ and $2(M-3)$ are the dimensions of the Grassmann varieties parametrizing the linear subspaces $L_b$. An easy analysis of the derivative shows that the function
\[ F(b) = b \cdot \frac{(M-1-b)(M-b)}{2} + b^2 + 1 \]
is everywhere increasing for $M \geq 5$, hence the minimum of $F(b)$ is $(M-2)(M-3) + 5$ if $2 \leq b \leq M-2$. Hence the overall minimum is
\[ \frac{M(M-3)}{2} + 3. \]
Following the same line of arguments we also obtain a lower bound for
\[ {\mathop{\rm codim}}_{\pi_{M,M-1}^{-1}(Q_{M-1})} \pi_{M,M-1}^{-1}(Q_{M-1}) \cap Z_M. \]
First note that it is not necessary to fix $q_1$ since linear terms do not occur. Hence $q_2, \ldots, q_M$ are polynomials in $X_1, \ldots, X_M$. Adapting the calculations above shows that the codimension is at least the minimum of the numbers
\[ b \cdot \frac{(M-b)(M+1-b)}{2} + b^2 + 1, 2 \leq b \leq M \]
and
\[ \frac{(M+1)(M-2)}{2} + 3, \]
that is $\frac{(M+1)(M-2)}{2} + 3$, arguing as before.

\vspace{0.3cm}
Finally, all these estimates imply that ${\mathop{\rm codim}}_{\cal F} \overline{{\cal F}\setminus{\cal F}_{\rm reg}}$ is bounded from below by the minimum of the numbers
\[ \binom{M}{i} - (M-1), 2 \leq i \leq M-2, \frac{M(M-3)}{2} + 3 - (M-1), \]
\[ \binom{M+1}{j} - (M-1) + M, 2 \leq j \leq M, \frac{(M+1)(M-2)}{2} + 3 - (M-1) + M, \]
that is
\[ \frac{M(M-3)}{2} + 3 - (M-1) = \frac{M(M-5)}{2} + 4 \]
for $M \geq 5$.\vspace{0.3cm}


\noindent {\bf 3. The $4n^2$-inequality.} Let us prove Proposition 1. We fix a mobile linear system $\Sigma$ on $V$ and a maximal singularity $E\subset V^+$ satisfying the Noether-Fano inequality $\mathop{\rm ord}_E\varphi^*\Sigma>na(E)$. We assume the $\mathop{\rm centre}B=\varphi(E)$ of $E$ on $V$ to be maximal, that is, $B$ is not contained in the centre of another maximal singularity of the system $\Sigma$. In other words, the pair $(V,\frac{1}{n}\Sigma)$ is canonical outside $B$ in a neighborhood of the generic point of $B$.\vspace{0.1cm}

Further, we assume that $B\subset\mathop{\rm Sing}V$ (otherwise the claim is well known), so that $\mathop{\rm codim}(B\subset V)\geq 4$. Let
$$
\begin{array}{cccl}
\varphi_{i,i-1}\colon & V_i &\to & V_{i-1}\\
& \cup & & \cup\\
& E_i & \to & B_{i-1}\\
\end{array}
$$
$i=1,\dots,K$, be the {\it resolution} of $E$, that is, $V_0=V$, $B_0=B$, $\varphi_{i,i-1}$ blows up $B_{i-1}=\mathop{\rm centre}(E,V_{i-1})$, $E_i=\varphi^{-1}_{i,i-1}(B_{i-1})$ the exceptional divisor, and, finally, the divisorial valuations, determined by $E$ and $E_K$, coincide.\vspace{0.1cm}

As explained in Sec. 4 below, for every $i=0,\dots,K-1$ there is a Zariski open subset $U_i\subset V_i$ such that $U_i\cap B_i\neq\emptyset$ is smooth and either $V_i$ is smooth along $U_i\cap B_i$, or every point $p\in U_i\cap B_i$ is a quadratic singularity of $V_i$ of rank at least 5. In particular, the quasi-projective varieties $\varphi^{-1}_{i,i-1}(U_{i-1})$, $i=1,\dots,K$, are factorial and the exceptional divisor
$$
E^*_i=E_i\cap\varphi^{-1}_{i,i-1}(U_{i-1})
$$
is either a projective bundle over $U_{i-1}\cap B_{i-1}$ (in the non-singular case) or a fibration into quadrics of rank $\geq 5$ over $U_{i-1}\cap B_{i-1}$ (in the singular case). We may assume that $U_i\subset\varphi^{-1}_{i,i-1}(U_{i-1})$ for $i=1,\dots,K-1$. The exceptional divisors $E^*_i$ are all irreducible.\vspace{0.1cm}

As usual, we break the sequence of blow ups into the {\it lower} $(1\leq i\leq L)$ and {\it upper} $(L+1\leq i\leq K)$ parts: $\mathop{\rm codim}B_{i-1}\geq 3$ if and only if $1\leq i\leq L$. It may occur that $L=K$ and the upper part is empty (see \cite{Pukh00a,Pukh98b,Pukh07b}). Set
$$
L_*=\mathop{\rm max}\{i=1,\dots,K\, |\, \mathop{\rm mult}\nolimits_{B_{i-1}}V_{i-1}=2\}.
$$
Obviously, $L_*\leq L$. Set also
$$
\delta_i=\mathop{\rm codim}B_{i-1}-2\quad \mbox{for}\quad 1\leq i\leq L_*
$$
and
$$
\delta_i=\mathop{\rm codim}B_{i-1}-1\quad \mbox{for} \quad L_*+1\leq i\leq K.
$$
We denote strict transforms on $V_i$ by adding the upper index $i$: say, $\Sigma^i$ means the strict transform of the system $\Sigma$ on $V_i$. Let $D\in\Sigma$ be a generic divisor. Obviously,
$$
D^i|_{U_i}=\varphi^*_{i,i-1}(D^{i-1}|_{U_{i-1}})-\nu_iE^*_i,
$$
where the integer coefficients $\nu_i=\frac12\mathop{\rm mult}_{B_{i-1}}\Sigma^{i-1}$ for $i=1,\dots,L^*$ and $\nu_i=\mathop{\rm mult}_{B_{i-1}}\Sigma^{i-1}$ for $i=L^*+1,\dots,K$.\vspace{0.1cm}

Now the Noether-Fano inequality takes the traditional form
\begin{equation}\label{19september2012.1}
\sum^K_{i=1}p_i\nu_i> n \left(\sum^K_{i=1}p_i\delta_i\right),
\end{equation}
where $p_i$ is the number of paths from the top vertex $E_K$ to the vertex $E_i$ in the oriented graph $\Gamma$ of the sequence of blow ups $\varphi_{i,i-1}$, see \cite{Pukh00a,Pukh98b,Pukh07b} for details.\vspace{0.1cm}

We may assume that $\nu_1<\sqrt{2}n$, otherwise for generic divisors $D_1,D_2\in\Sigma$ we have
$$
\mathop{\rm mult}\nolimits_B(D_1\circ D_2)\geq 2\nu^2_1>4n^2
$$
and the $4n^2$-inequality is shown. We do not use the following claim, but nevertheless it is worth mentioning.\vspace{0.1cm}

{\bf Lemma 1.} {\it The inequality $\nu_1>n$ holds.}\vspace{0.1cm}

{\bf Proof.} Taking a point $p\in B$ of general position and a generic complete intersection 3-germ $Y\ni p$, we reduce to the case of a non log canonical singularity centered at a non-degenerate quadratic point, when the claim is well known, see \cite{Corti00,Pukh09c}. Q.E.D.\vspace{0.1cm}

Obviously, the multiplicities $\nu_i$ satisfy the inequalities
\begin{equation}\label{19september2012.2}
\nu_1\geq\dots\geq\nu_{L^*}
\end{equation}
and, if $K\geq L^*+1$, then
\begin{equation}\label{19september2012.3}
2\nu_{L^*}\geq\nu_{L^*+1}\geq\dots\geq\nu_K.
\end{equation}
Now let $Z=(D_1\circ D_2)$ be the self-intersection of the mobile system $\Sigma$ and set $m_i=\mathop{\rm mult}_{B_{i-1}}Z^{i-1}$ for $1\leq i\leq L$. Applying the technique of counting multiplicities in word for word the same way as in \cite{Pukh00a,Pukh98b,Pukh07b}, we obtain the estimate
$$
\sum^L_{i=1}p_im_i\geq 2\sum^{L_*}_{i=1}p_i\nu^2_i+\sum^K_{i=L_*+1}p_i\nu^2_i.
$$
Denote the right hand side of this inequality by $q(\nu_1,\dots,\nu_K)$. We see that
$$
\sum^L_{i=1}p_im_i>\mu,
$$
where $\mu$ is the minimum of the positive definite quadratic form $q(\nu_*)$ on the compact convex polytope $\Delta$ defined on the hyperplane
$$
\Pi=\left\{\sum^K_{i=1}p_i\nu_i=n\left(\sum^K_{i=1}p_i\delta_i\right)\right\}
$$
by the inequalities (\ref{19september2012.2},\ref{19september2012.3}). Let us  estimate $\mu$.\vspace{0.1cm}

We use the standard optimization technique in two steps. First, we minimize $q|_{\Pi}$ separately for the two groups of variables
$$
\nu_1,\dots,\nu_{L_*}\quad \mbox{and}\quad \nu_{L_*+1},\dots,\nu_K.
$$
Easy computations show that the minimum is attained for
$$\nu_1=\dots=\nu_{L_*}=\theta_1\quad\mbox{and}\quad\nu_{L_*+1}=\dots=\nu_K=\theta_2
,
$$
satisfying the inequality $2\theta_1\geq\theta_2$. Putting
$$
\Sigma_*=\sum^{L_*}_{i=1}p_i\quad \mbox{and}\quad \Sigma^*=\sum^K_{i=L_*+1}p_i,
$$we get the extremal problem
$$
\bar{q}(\theta_1,\theta_2)=2\Sigma_*\theta^2_1+\Sigma^*\theta^2_2\to\mathop{\rm min}
$$
on the ray, defined by the inequality $2\theta_1\geq\theta_2$ on the line
$$
\Lambda=\left\{\Sigma_*\theta_1+\Sigma^*\theta_2=n\sum^K_{i=1}p_i\delta_i\right\}.
$$
Now we make the second step, minimizing $\bar{q}|_{\Lambda}$. The minimum is attained for $\theta_1=\theta$, $\theta_2=2\theta$ (so that the condition $2\theta_1\geq\theta_2$ is satisfied and for that reason can be ignored), where $\theta$ is obtained from the equation of the line $\Lambda$:
$$
\theta=\frac{n}{\Sigma_*+2\Sigma^*}\sum^K_{i=1}p_i\delta_i.
$$
Now set
$$
\Sigma_l=\sum^L_{i=1}p_i,\quad \Sigma^*_l=\sum^L_{i=L_*+1}p_i,\quad \Sigma_u=\sum^K_{i=L+1}p_i
$$
(if $L\geq L_*+1$; otherwise set $\Sigma^*_l=0$). Obviously, the relations
\begin{equation}\label{19september2012.4}
\Sigma_l=\Sigma_*+\Sigma^*_l\quad \mbox{and}\quad \Sigma^*=\Sigma^*_l+\Sigma_u
\end{equation}
hold. Recall that, due to our assumptions on the singularities of $V_i$ we have $\delta_i\geq 2$ for $i\leq L$. Therefore,
$$
\theta\geq \frac{2\Sigma_l+\Sigma_u}{\Sigma_*+2\Sigma^*}n
$$
and so
$$
\mu\geq 2\frac{(2\Sigma_l+\Sigma_u)^2}{\Sigma_*+2\Sigma^*}n^2.
$$
Since
$$
\Sigma_l\mathop{\rm mult}\nolimits_BZ\geq\sum^L_{i=1}p_im_i,
$$we finally obtain the estimate
$$
\mathop{\rm mult}\nolimits_BZ>2\frac{(2\Sigma_l+\Sigma_u)^2}{\Sigma_l(\Sigma_*+2\Sigma^*)}n^2.
$$
Therefore, the $4n^2$-inequality follows from the estimate
$$
(2\Sigma_l+\Sigma_u)^2\geq 2\Sigma_l(\Sigma_*+2\Sigma^*).
$$
Replacing in the right hand side $\Sigma_*+2\Sigma^*$ by
$$
\Sigma_*+2(\Sigma^*_l+\Sigma_u)=\Sigma_l+\Sigma^*_l+2\Sigma_u,
$$
we bring the required estimate to the following form:
$$
2\Sigma^2_l+\Sigma^2_u\geq 2\Sigma_l\Sigma^*_l,
$$
which is an obvious inequality. Proof of Proposition 1 is now complete. Q.E.D.\vspace{0.3cm}


\noindent {\bf 4. Stability of the quadratic singularities under blow ups.} We start with the following essential\vspace{0.1cm}

{\bf Definition 1.}
Let $X \subset Y$ be a subvariety of codimension $1$ in a smooth quasi-projective complex variety $Y$ of dimension $n$. A point $P \in X$ is called a quadratic point of rank $r$ if there are analytic coordinates $z = (z_1, \ldots, z_n)$ of $Y$ around $P$ and a quadratic form $q_2(z)$ of rank $r$ such that the germ of $X$ in $P$ is given by
\[ (P \in X) \cong \left\{ q_2(z) + \mathrm{terms\ of\ higher\ degree} = 0 \right\} \subset Y. \]

\vspace{0.3cm}
{\bf Theorem 4.}
Let $X \subset Y$ be a subvariety of codimension $1$ in a smooth quasi-projective complex variety $Y$ of dimension $n$, with at most quadratic points of rank $\geq r$ as singularities. Let $B \subset X$ be an irreducible subvariety. Then there exists an  open set $U \subset Y$ such that
\begin{itemize}
\item[(i)] $B \cap U$ is smooth, and
\item[(ii)] the blow up $\widetilde{X}_U$ of $X \cap U$ along $B \cap U$ has at most quadratic points of rank $\geq r$ as singularities.
\end{itemize}

\begin{proof}
The statement is obvious if $B \not\subset \mathrm{Sing}(X)$. So we assume from now on that $B \subset \mathrm{Sing}(X)$.

By restricting to a Zariski-open subset of $Y$ we may assume that $B \subset \mathrm{Sing}(X)$ is a smooth subvariety. By assumption there exist analytic coordinates $z = (z_1, \ldots, z_n)$ around each point $P \in B \subset Y$ such that the germ
\[ (P \in X) \cong \left\{ f(z) = z_1^2 + \ldots + z_r^2 + \mathrm{terms\ of\ higher\ degree} = 0 \right\} \subset Y. \]
Then the singular locus $\mathrm{Sing}(X)$ is contained in the vanishing locus of the partial derivatives of this equation, hence in
\[ \left\{ \frac{\partial f}{\partial z_1} = \cdots = \frac{\partial f}{\partial z_r} = 0 \right\}. \]
Since
\[ \frac{\partial f}{\partial z_i} = 2z_i + \mathrm{terms\ of\ higher\ degree}, 1 \leq i \leq r, \]
setting $z^\prime_1 := \frac{1}{2}\frac{\partial f}{\partial z_1}, \ldots, z^\prime_r := \frac{1}{2}\frac{\partial f}{\partial z_r}, z^\prime_i := z_i$ for $r+1 \leq i \leq n$ yields new analytic coordinates
\[ z^\prime_1, \ldots, z^\prime_r, z^\prime_{r+1}, \ldots, z^\prime_n \]
of $Y$ around $P$. In these new coordinates the defining equation of $X$ still is of the form
\[ {z^\prime_1}^2 + \ldots + {z^\prime_r}^2 + \mathrm{terms\ of\ higher\ degree} = 0, \]
and $B \subset \left\{ z^\prime_1 = \ldots = z^\prime_r = 0 \right\}$. Perhaps after a further coordinate change we can even assume that
\[ B = \left\{ z^\prime_1 = \ldots = z^\prime_k = 0 \right\}, k \geq r. \]

\vspace{0.2cm}

\noindent \textit{Claim.} $(P \in X) \cong \left\{ {z^\prime_1}^2 + \ldots + {z^\prime_r}^2 + f_{\geq 3} = 0 \right\}$ where $f_{\geq 3}$ consists of terms of degree $\geq 3$ and is an element of $(z^\prime_1, \ldots, z^\prime_k)^2$.

\vspace{0.1cm}

\noindent \textit{Proof of Claim.} $B \subset \mathrm{Sing}(X)$ must be contained in $\left\{ \frac{\partial f_{\geq 3}}{\partial z^\prime_j} = 0 \right\}$, hence $\frac{\partial f_{\geq 3}}{\partial z^\prime_j} \in (z^\prime_1, \ldots, z^\prime_k)$ for all $k+1 \leq j \leq n$. This is only possible if $f_{\geq 3} \in (z^\prime_1, \ldots, z^\prime_k)$. Write $f_{\geq 3} = z^\prime_1 f^\prime_1 + \ldots + z^\prime_k f^\prime_k$. Then as before $\frac{\partial f_{\geq 3}}{\partial z^\prime_i} = f^\prime_i + \sum_{1 \leq j \leq k, j \neq i} z^\prime_j \frac{\partial f^\prime_j}{\partial z^\prime_i} \in (z^\prime_1, \ldots, z^\prime_k)$ for all $1 \leq i \leq k$. But this is only possible if $f^\prime_i \in (z^\prime_1, \ldots, z^\prime_k)$ for all $1 \leq i \leq k$. $\hfill$ $\Box$

\vspace{0.2cm}

\noindent Using the coordinates $z^\prime_1, \ldots, z^\prime_n$ we can cover the blow up of $Y$ along $B$ over $P \in Y$ by $k$ charts with coordinates
\[ t_1^{(i)}, \ldots, z_i, \ldots, t_k^{(i)}, z_{k+1}, \ldots, z_n, 1 \leq i \leq k, \]
where $z^\prime_j = t_j^{(i)} z_i$ for $1 \leq j \leq k$, $j \neq i$, $z^\prime_i = z_i$ and $z^\prime_l = z_l$ for $k+1 \leq l \leq n$.
To prove the theorem we only need to check in each chart that along the fiber of the exceptional divisor over $P \in B$ there are at most quadratic points of rank $\geq r$ as singularities. We distinguish several cases:

\vspace{0.2cm}

\noindent \textit{Case 1.} $1 \leq i \leq r$, say $i = 1$.

\noindent Then the strict transform of $X$ is given by the equation
\[ 1 + (t_2^{(1)})^2 + \cdots + (t_r^{(1)})^2 + z_1 \cdot F + Q(t_2^{(1)}, \ldots, t_k^{(1)}) \cdot G = 0, \]
where $Q$ is a quadratic polynomial in $t_2^{(1)}, \ldots, t_k^{(1)}$ and $G \in (z_{k+1}, \ldots, z_n)$.
On the fiber of the exceptional divisor over $P$, $\left\{ z_1 = z_{k+1} = \ldots = z_n = 0 \right\}$, the gradient of this function can only vanish when $t_2^{(1)} = \ldots = t_r^{(1)} = 0$. But this locus does not intersect the strict transform, hence in this chart the strict transform is smooth along the fiber of the exceptional divisor over $P$.

\vspace{0.2cm}

\noindent \textit{Case 2.} $r+1 \leq i \leq k$, say $i = k$.

\noindent Then the strict transform of $X$ is given by the equation
\[ (t_1^{(k)})^2 + \cdots + (t_r^{(k)})^2 + z_k \cdot F + Q(t_1^{(k)}, \ldots, t_{k-1}^{(k)}) \cdot G = 0, \]
$Q$ and $G$ as above. On the fiber of the exceptional divisor over $P$, $\left\{ z_k = z_{k+1} = \ldots = z_n = 0 \right\}$, the gradient of this function can only vanish when $t_1^{(k)} = \ldots = t_r^{(k)} = 0$. We first discuss the origin in these coordinates,
\[ (0, \ldots, 0) \in \{ t_1^{(k)} = \ldots = t_r^{(k)} = z_k = z_{k+1} = \ldots = z_n = 0 \}. \]
If $F$ has a constant term then the strict transform of $X$ is smooth in $(0, \ldots, 0)$.

\noindent If $F$ has no constant terms but contains linear terms then the rank of the quadratic term in the defining equation is still $\geq r$ because we only add quadratic monomials containing $z_k$ to $(t_1^{(k)})^2 + \cdots + (t_r^{(k)})^2$. Hence $(0, \ldots, 0)$ is a quadratic point of rank $\geq r$.

\noindent Finally, if $F$ is of degree $\geq 2$ the quadratic term in the defining equation is $(t_1^{(k)})^2 + \cdots + (t_r^{(k)})^2$. Hence $(0, \ldots, 0)$ is a quadratic point of rank $r$.

\noindent The affine coordinate change to
\[ t_1^{(k)}, \ldots, t_r^{(k)}, t_{r+1}^{(k)} - a_{r+1}, \ldots, t_{k-1}^{(k)} - a_{k-1}, z_k, z_{k+1}, \ldots, z_n \]
leads to a defining equation of the strict transform around the point
\[ (0, \ldots, 0, a_{r+1}, \ldots, a_{k-1}, 0, 0, \ldots, 0) \in \{ t_1^{(k)} = \ldots = t_r^{(k)} = z_k = z_{k+1} = \ldots = z_n = 0 \} \]
in one of the forms already discussed. Consequently, in this chart all points in the strict transform of $X$ also lying on the fiber of the exceptional divisor over $P$ are smooth or quadratic points of rank $\geq r$.
\end{proof}

\vspace{0.3cm}
{\bf Remark 2.}
Note that $\widetilde{X}_U$ is again a subvariety of codimension $1$ in the smooth quasi-projective blow up of $U$ along $B \cap U$. The universal property of blow ups~\cite[Prop.II.7.14]{Hart:AG} and the calculations in the proof above tell us that the exceptional locus $E_U \subset \widetilde{X}_U$ is a Cartier divisor on $\widetilde{X}_U$ such that the morphism $E_U \rightarrow B \cap U$ is a fibration into quadrics of rank $\geq r$ in a $\mathbb{P}^{\mathrm{codim}_Y B}$-bundle.


\end{document}